\documentclass[alpha]{article}

\usepackage{amssymb}
\usepackage{amsmath,amsfonts}
\usepackage{theorem}
\usepackage[dvips]{graphicx}
\usepackage{graphics}
\usepackage{epsfig}
\usepackage{alltt}
\usepackage[english]{babel}

\newtheorem{theorem}{Theorem}[section]
\newtheorem{definition}[theorem]{Definition}

\newtheorem{lemma}[theorem]{Lemma}
{\theorembodyfont{\rmfamily}} %sense italica
{\theorembodyfont{\rmfamily}} %sense italica

\def\card#1{\left|#1\right|}

\usepackage{pst-all}
\usepackage{pstcol}
\usepackage{subfigure}
\psset{xunit=.2cm, yunit=.2cm}

\title{A Fibonacci sequence for linear structures\\
with two types of components\footnote{All
the results contained in this file are included in a paper submitted to
Annals of Operations Research in October, 2008 on ocasion of the 
Conference on Applied Mathematical Programming and Modelling,
that held in Bratislava in May, 2008.}}

\author{J. Freixas\thanks{Universitat Polit\`ecnica de Catalunya.
DMA3 and EPSEM. E-08240 Manresa, Spain. \texttt{josep.freixas@upc.edu}}
        \and
        X. Molinero\thanks{Universitat Polit\`ecnica de Catalunya.
LSI and EPSEM. E-08240 Manresa, Spain. \texttt{molinero@lsi.upc.edu}}
        \and
        S. Roura\thanks{Universitat Polit\`ecnica de Catalunya.
LSI. E-08034 Barcelona, Catalonia, Spain. \texttt{roura@lsi.upc.edu}}
}

\date{\today}

\begin{document}

\maketitle

%\thispagestyle{empty}
%\newpage

\begin{abstract}
We investigate binary voting systems with two types of voters and a
hierarchy among the members in each type, so that members in one
class have more influence or importance than members in the other
class. 
The purpose of this paper is to count,
up to isomorphism, the number of these voting systems for an
arbitrary number of voters.
We obtain a closed formula for the number of these systems,
this formula follows a Fibonacci sequence with a smooth polynomial
variation on the number of voters. 
\newline


\end{abstract}

\bigskip
\section{Introduction}
\label{sec:intro} We consider voting systems in which each player
casts a ``yes" or ``no" vote, and the outcome is a collective ``yes"
or ``no". These voting systems, known in the literature as simple
games, can be very complicated. Their specialization to
\emph{symmetric} simple games, however, are simple indeed; each such
game corresponds to the qualified majority rule, in which an issue is
passed if and only if the number of voters in favour meets or exceeds
some \emph{threshold} or \emph{quota}. We refer, as in
\cite{FrZw08GEB}, to this result as \emph{May's Theorem for Simple
Games} ``with bias"\footnote{Because May's original result (see
\cite{May52} or \cite{Tay95}) considers only anonymous voting
systems (for two alternatives) that are \emph{neutral}: there is no
built-in bias towards ``yes" or ``no" outcomes, so that they are
treated symmetrically.}. Thus, May's Theorem with bias may be
stated as follows: for each positive integer $n$, there is, up to
isomorphism, a unique simple game that is anonymous or symmetric
(i.e., voters play an equivalent r\^{o}le in the game).
At least two facts are relevant of this result:
\begin{enumerate}
\item all symmetric games are weighted,\footnote{Loosely speaking,
a game is weighted if it can be assigned a quota and a weight for
each voter, so that winning coalitions are those with the sum of the
weights of their members greater than or equal to the quota.}
\item the function $S$ on the number of voters $n$ that counts all these
games is, up to isomorphism, the identity, i.e., $S(n)=n$.
\end{enumerate}

There seems to be only one natural direction in which to extend
symmetric games\footnote{Symmetric games are also called
$q$--out--of--$n$ games, understanding that $n$ is the number of
voters, a weight of $1$ is assigned to each voter, and the quota $q$
is a fixed integer between $1$ and $n$.} within
simple games. To accommodate this new class of voting systems, we
make two changes to the symmetric simple games: we allow
\emph{two} classes of symmetric voters
instead of one, and we consider that voters in
one class are more influential or important than voters in the other
class. The goal of this paper is to analyze the joint effect of
these two changes.
What voting systems are possible?
How many of them are these, up to isomorphism, for a fixed number of voters?

As we shall see in this paper, the obtained results are paradoxically
opposed to those described above for symmetric games, which
demonstrates that the complexity of these close voting systems is
considerably higher than that of symmetric games. Indeed, we will prove
that the number of them is $F(n+6)-(n^2 + 4n +8)$, where $F(n)$ are
the Fibonacci numbers\footnote{The Fibonacci numbers are
defined by the following recurrence relation: $F(0)=0$,
$F(1)=1$, and $F(n)=F(n-1)+F(n-2)$ for all $n>1$.}, being $n$ the
number of voters. Hence, the number of these games is asymptotically
exponential, which contrasts with the linear behavior of symmetric
games.

The paper is organized as follows. Basic definitions and preliminary
results are included in Section~\ref{sec:preliminaries}.
Section~\ref{sec:counting} contains the main result of the paper, which is
devoted to count the number of voting systems in terms of the number
of voters, resulting a closed formula with asymptotic exponential
behavior.

\section{Preliminaries}
\label{sec:preliminaries}
Simple games can be viewed as models of voting systems in which a
single alternative, such as a bill or an amendment, is pitted
against the status quo.

\begin{definition}
\label{DSG} A \emph{simple game} is a pair $(N, W)$ in which $N =
\{1,2,\dots,n\}$, and $W$ is a collection of subsets of $N$ that
satisfies $N \in W$, $\emptyset \notin W$ and the monotonicity
property: if $S\in{}W$ and $S\subseteq{}T\subseteq{}N$, then
$T\in{}W$.
\end{definition}

Any set of voters is called a \emph{coalition}, and the set $N$ is
called the \emph{grand coalition}. Members of $N$ are called
\textit{players} or \textit{voters}, and the subsets of $N$ that are
in $W$ are called \emph{winning coalitions}. The subfamily of
\emph{minimal} winning coalitions $W^m = \{S\in W: T\subset
S\Rightarrow T\notin W\}$ determines the game. The subsets of $N$
that are not in $W$ are called \emph{losing coalitions}. The
subfamily of \emph{maximal} losing coalitions is $L^M = \{S\in L:
S\subset T\Rightarrow T\in W\}$. A voter $i$ is \emph{null} in
$(N,W)$ if $i \notin S$ for all $S \in W^m$. Thus, for a non--null
voter $i$ there is at least a coalition $S$ such that $S\in W$, $i
\in S$ and $S\setminus \{i\} \in L$. Real--world examples of simple
games are given by Taylor~\cite{Tay95,TaZw99}.

The ``desirability" relation defined on the set of voters represents
a way to make precise the idea that a particular voting system may
give one voter more influence than another. Isbell already used it
in \cite{Isb58}.

\begin{definition}
\label{Ddes} Let $(N,W)$ be a simple game.

\begin{itemize}
\item[$(i)$] Player $i$ is \emph{more desirable} than $j$
($i \succsim  j$, in short) in $(N,W)$ iff
$$ S\cup\{j\}\in W\ \Rightarrow   \ S\cup\{i\}\in W,\qquad
   \text{for all $S\subseteq N \setminus \{i,j\}$}.$$

\item[$(ii)$] Players $i$ and $j$ are \emph{equally desirable}
($i \approx j$, in short) in $(N,W)$ iff
$$ S\cup\{i\}\in W \Leftrightarrow\ S\cup\{j\}\in W,\qquad
   \text{for all $S \subseteq N\setminus \{i,j\}$}.$$

\item[$(iii)$] Player $i$ is \emph{strictly more desirable}
than player $j$ ($i \succ j$, in short) in $(N,W)$ iff $i$ is more
desirable than $j$, but $i$ and $j$ are not equally desirable.
\end{itemize}
\end{definition}

\begin{definition}
\label{Dcomp} A simple game $(N,W)$ is \emph{complete} or
\emph{linear} if the desirability relation is a complete
preordering.
\end{definition}

In the field of Boolean algebra, complete games correspond to
$2$-monotonic positive Boolean functions, which were already
considered in \cite{Hu65}. The problem of identifying this type of
functions by using polynomial-time recognition have been treated in
\cite{BHIK91, BHIK97}. In a complete simple game we may decompose
$N$ in a collection of subsets, called classes, $N_1 > N_2
> \dots > N_t$, forming a partition of $N$, and understanding that if
$i \in N_p$ and $j \in N_q$ then: $p=q$ if and only if $i \approx j$;
and $p<q$ if and only if $i \succ j$.

Now we are going to define the $\delta$-ordering, introduced in
\cite{CaFr96} for an arbitrary number of types of voters.

\begin{definition} \label{Dprofiles}
Let $n$ be the number of players of a complete simple game $(N,W)$
with two types of players $N_1 > N_2$. Let $n_1 = \card{N_1}$ and
let $n_2=\card{N_2}$, where $(n_1,n_2)\in\mathbb{N}\times\mathbb{N}$
with $n_1+n_2=n$. Then the \emph{rectangle of $(n_1+1)\times(n_2+1)$ profiles}
for $(N,W)$ is:
$$I_{n_1} \times I_{n_2} = \{(m_1,m_2) \in (\mathbb{N} \cup \{ 0
\})\times (\mathbb{N} \cup \{ 0 \}) \, : \, m_1 \leq n_1, m_2\leq
n_2 \}.$$
In $I_{n_1} \times I_{n_2}$, the \emph{$\delta$-ordering}
given by the comparison of partial sums is:
$$ (p_1,p_2) \, \delta \, (m_1,m_2) \quad \emph{if and only if}
\quad p_1 \geq m_1 \; \; \text{and} \; \; p_1+p_2 \geq m_1+m_2.$$
\end{definition}

It is not difficult to check that the pair $(I_{n_1} \times
I_{n_2}, \delta)$ is a distributive lattice that possesses a maximum
element $(n_1, n_2)$ and a minimum element $(0,0)$. The profiles
in $I_{n_1} \times I_{n_2}$ can be completely ordered by the
\emph{lexicographical ordering}: profile $(p_1,p_2)$ is
lexicographically greater than $(m_1,m_2)$ if either $p_1>m_1$, or
$p_1=m_1$ with $p_2 > m_2$.

\begin{definition} \label{Disom}
Two simple games $(N, W)$ and $(N', W')$ are said to be
\emph{isomorphic} if there is a bijective map $f:N\rightarrow N'$
such that $S \in W$ if and only if $f(S) \in W'$; $f$ is called an
isomorphism of simple games.
\end{definition}

The following known result has three parts. The first part shows how
to associate a vector $(n_1,n_2)$ and a matrix $\mathcal M$  to a
complete simple game $(N,W)$, and describes the restrictions that
these parameters need to fulfill. The second part establishes that
every pair of isomorphic complete simple games $(N,W)$ and $(N',W')$
corresponds to the same associated vector $(n_1,n_2)$ and matrix $\mathcal M$
(uniqueness). The third part shows that a vector $(n_1,n_2)$ and a
matrix $\mathcal M$ fulfilling the conditions in Part A always correspond
to a complete simple game $(N,W)$ (existence).

\begin{theorem} (Carreras and Freixas' Theorems 4.1 and 4.2  in \cite{CaFr96} for $2$ types of voters)
\label{TCaFr2types}
\begin{description}
       \item[Part A] Let $(N,W)$ be a complete simple
       game with two nonempty classes $N_1>N_2$, and let $ (n_1,n_2)$ be
       the vector defined by their cardinalities. For each coalition $S \in
       W$, consider the node or profile $(s_1,s_2) \in I_{n_1} \times I_{n_2}$
       with components $s_k = \card{S \cap N_k}$ ($k=1,2$). Let $\mathcal
       M$ be the matrix
       $$\mathcal M=\left(%
\begin{array}{cc}
  m_{1,1} & m_{1,2} \\
\vdots & \vdots \\
  m_{r,1} & m_{r,2} \\
\end{array}%
\right)$$ whose $r$ rows are the nodes corresponding to
       winning coalitions which are minimal in the
       $\delta$-ordering. Matrix $\mathcal M$
       satisfies the two conditions below:
\begin{enumerate}
      \item $m_{i,1},m_{i,2} \in \mathbb{N} \cup \{0\}$,
            $0 \leq m_{i,1} \leq n_{1}$ and
            $0 \leq m_{i,2} \leq n_{2}$
            for all $1 \leq i \leq r$;
      \item if $r=1$, then $m_{1,1}>0$ and $m_{1,2}<n_2$;\par
            if $r\ge2$, then $m_{i,1}>m_{j,1}$ and
                             $m_{i,1}+m_{i,2}<m_{j,1}+m_{j,2}$
                             for all
                             $1\le{}i<j\le{}r$.\footnote{The
                             lexicographic ordering chosen
                             guarantees uniqueness under
                             permutation of rows.
                             This lexicographic ordering is a plausible
                             choice which could be replaced
                             for other alternative criteria.}
\end{enumerate}
     \item[Part B] (Uniqueness) Two complete simple games with two types of voters
     $(N, W)$ and $(N', \mathcal
     W')$ are isomorphic if and only if $(n_1,n_2) = (n'_1,n'_2)$
     and $\mathcal M = \mathcal M'$.
     \item[Part C] (Existence) Given a vector $(n_1,n_2)$ and a
     matrix $\mathcal M$ satisfying the conditions of Part A, there
     exists a complete simple game $(N,W)$ with two types of voters
     associated to $(n_1,n_2)$ and $\mathcal M$.
\end{description}
\end{theorem}

For example, to illustrate how $(N,W)$ is obtained,
let $(n_1,n_2)=(2,3)$ and $\mathcal M = \left(%
\begin{array}{cc}
  2 & 0 \\
  0 & 3 \\
\end{array}%
\right).$ The set of winning profiles in $I_2 \times I_3$ is
$\{(0,3), (1,2), (1,3), (2,0), (2,1), (2,3) \}$ because each of
these profiles either $\delta$-dominates $(2,0)$ or $(0,3)$; and the
set of minimal winning profiles is $\{ (0,3), (1,2), (2,0) \}$,
because the other winning profiles can be obtained from one of these
three profiles by simply adding some element in any of its
components. If we take $N_1=\{1,2\}$ and $N_2=\{3,4,5\}$, then $N =
N_1 \cup N_2$, and
$$W^m = \{ \{1,2\}, \{1,3,4\}, \{1,3,5\}, \{1,4,5\}, \{2,3,4\},
\{2,3,5\}, \{2,4,5\}, \{3,4,5\} \},$$ where the first coalition
corresponds to profile $(2,0)$, the last coalition to profile
$(0,3)$ and the remaining coalitions to profile $(1,2)$.

Theorem \ref{TCaFr2types} is a \emph{parametrization} theorem
because it allows one to enumerate all complete games up to
isomorphism by listing the possible values of certain invariants.
We will see in next section that such enumeration can be done for
every number of voters.

\section{Counting complete games with two types of voters: Fibonacci numbers}
\label{sec:counting}
We establish a relation between the number of non--isomorphic
complete simple games with $n$ voters of two different types and
the Fibonacci numbers.

\begin{definition}
\label{def:H(n)}
Let $H(n)$ denote the number of non--isomorphic complete simple games
with $n$ voters of two different types.
\end{definition}

Note that $H(n)$ is the number of matrices
\[
\left(
 \begin{matrix}
  m_{1,1}\hspace{.25cm} & m_{1,2}\\
  m_{2,1}\hspace{.25cm} & m_{2,2}\\
  \vdots\hspace{.25cm}  & \vdots\\
  m_{r,1}\hspace{.25cm} & m_{r,2}
 \end{matrix}
\right)
\]
such that there exists $(n_1,n_2)\in\mathbb{N} \times \mathbb{N}$,
with $n=n_1+n_2$, verifying Properties $(1)-(2)$ of
Theorem~\ref{TCaFr2types}-(A).

\begin{theorem}
\label{th:H(n)}
$$ H(n) = F(n+6) - (n^2+4n+8), $$
where $F(n)$ is the $n$-th Fibonacci number.
\end{theorem}

The proof of this theorem is a consequence of some
additional definitions and lemmas.

\begin{definition}
\label{def:N(a,s,b)}
Let $N(a,s,b)$ be the number of matrices with non-negative integer entries
and an arbitrary number of rows $r\ge1$
\[
\left(
 \begin{matrix}
  m_{1,1}\hspace{.25cm} & m_{1,2}\\
  m_{2,1}\hspace{.25cm} & m_{2,2}\\
  \vdots\hspace{.25cm}  & \vdots\\
  m_{r,1}\hspace{.25cm} & m_{r,2}
 \end{matrix}
\right)
\]
such that
\begin{enumerate}
  \item $m_{1,1}=a$;
  \item $m_{1,1} + m_{1,2} = s$;
  \item $m_{r,2}=b$;
  \item if $r\ge2$, then
        $m_{i,1} > m_{j,1}$ and $m_{i,1}+m_{i,2} < m_{j,1}+m_{j,2}$
        for all $i$ and $j$ with $1\le{}i<j\le{}r$.
\end{enumerate}
\end{definition}

\begin{lemma}
\label{lem:N(a,s,b)}
Let $N(a,s,b)$ be as defined above, then
\[
 N(a,s,b) =
 \begin{cases}
  0, & \text{if $s-a>b$};\\
  1, & \text{if $s-a=b$};\\
  \sum\limits_{k=0}^{a-1}{{b+a-2-s}\choose{k}}, & \text{if $s-a<b$}.
 \end{cases}
\]
\end{lemma}

\noindent \emph{Proof:}
We consider three cases depending on the
relation between $a$, $s$ and $b$.

\noindent{}\emph{Case 1: $s-a>b$.}\par
By definition, it is clear that $N(a,s,b)$ is equal to $0$ whenever $s-a>b$.

\noindent{}\emph{Case 2: $s-a=b$.}\par
Note that $s-a=b$ implies $r=1$.
In fact, we just have one matrix with one row
\[ (a\ \ s-a)\equiv(a\ \ b). \]
Thus, $N(a,s,s-a)=1$.

\noindent{}\emph{Case 3: $s-a<b$.}\par
It is clear that $s-a<b$ implies $r>1$, and
$N(a,s,b)$ is equal to the number of matrices with $r-1$ rows
\[
\left(
 \begin{matrix}
  m_{2,1}\hspace{.25cm} & m_{2,2}\\
  \vdots\hspace{.25cm}  & \vdots\\
  m_{r,1}\hspace{.25cm} & b
 \end{matrix}
\right)
\]
such that verify (from Definition~\ref{def:N(a,s,b)}):
\begin{enumerate}
 \item $0\le{}m_{2,1}<a$;
 \item $m_{2,1}+m_{2,2}>s$;
 \item if $r\ge3$, then $m_{i,1} > m_{j,1}$ and
                        $m_{i,1}+m_{i,2} < m_{j,1}+m_{j,2}$
                        for all $i$ and $j$ with $2\le{}i<j\le{}r$.
\end{enumerate}

That is,
\[
 N(a,s,b) = \sum\limits_{i=0}^{a-1}
             \sum\limits_{j>s}
              N(i,j,b).
\]
On the other hand, since $N(i,j,b)=0$ if $j>i+b$ (\emph{Case 1}),
we have that
\[
 N(a,s,b) = \sum\limits_{i=0}^{a-1}
             \sum\limits_{j=s+1}^{i+b}
              N(i,j,b).
\]

Now, by \emph{mathematical induction} we are going to prove the equality
\[
 \sum\limits_{i=0}^{a-1}\sum\limits_{j=s+1}^{i+b} N(i,j,b) =
 \sum\limits_{k=0}^{a-1} {{b+a-2-s}\choose{k}}.
\]

\noindent{}By \emph{induction hypothesis} [\textit{i.h.}], it
is clear that
\[
{\allowdisplaybreaks
\begin{array}{rcl}
  \sum\limits_{i=0}^{a-1}
   \sum\limits_{j=s+1}^{i+b} N(i,j,b)
&=&
  \sum\limits_{i=0}^{a-1}
   \sum\limits_{j=s+1}^{b+i-1} N(i,j,b) +
  f(a,s,b) \\
&{[\textit{i.h.}]\atop{}=}&
  \sum\limits_{i=0}^{a-1}
   \sum\limits_{j=s+1}^{b+i-1}
    \sum\limits_{k=0}^{i-1} {{b+i-2-j}\choose{k}} +
  f(a,s,b)
\end{array}
}\]
\noindent{}where we define
\begin{equation}
  f(a,s,b) := \sum\limits_{i=0}^{a-1} N(i, b+i, b)
              \quad\text{if and only if $(s + 1 \le b + i)$}.\notag
\end{equation}

\noindent{}Note that $s+1\le{}b+i \iff s+1-b\le{}i$;
thus, it is clear that (by definition):
\begin{equation}
f(a,s,b) =
\begin{cases}
 a,       & \text{if $(s<b)     \iff (s+1-b\le0)$} \\
          & \label{eq:f(a,s,b)}\\
 a+b-1-s, & \text{if $(s\ge{}b) \iff (s+1-b>{}0)$}
\end{cases}
\end{equation}

From now on we will use some elementary combinatorial equalities
like~\cite{Cam94}
\[
 \sum\limits_{\alpha\le\beta}{\alpha\choose\gamma} =
 {{\beta+1}\choose{\gamma+1}},
 \hspace{.5cm}\text{for all $\gamma\ge{}0$.}
\]

First, we consider the following equalities
\[
{\allowdisplaybreaks
 \begin{array}{rl}
  \sum\limits_{j=s+1}^{b+i-1}
   \sum\limits_{k=0}^{i-1}
    {{b+i-2-j}\choose{k}} =&
  \sum\limits_{j\ge{}s+1}
   \sum\limits_{k\le{}i-1}
    {{b+i-2-j}\choose{k}}
 =
  \sum\limits_{k\le{}i-1}
   \sum\limits_{j\ge{}s+1}
    {{b+i-2-j}\choose{k}}\\
 =&
  \sum\limits_{k\le{}i-1}
   \sum\limits_{z\le{}b+i-3-s}
    {{z}\choose{k}}
 =
  \sum\limits_{k=0}^{i-1}
   {{b+i-2-s}\choose{k+1}}
 \end{array}
}\]

Second, reorganizing the previous results we have:
\[
{\allowdisplaybreaks
 \begin{array}{rl}
  \sum\limits_{i=0}^{a-1}
   \sum\limits_{k=0}^{i-1}
    {{b+i-2-s}\choose{k+1}}
 =&
  \sum\limits_{k=0}^{a-2}
   \sum\limits_{i=k+1}^{a-1}
    {{b+i-2-s}\choose{k+1}}
 =
  \sum\limits_{k=0}^{a-2}
   \sum\limits_{y=b+k-1-s}^{b+a-3-s}
    {{y}\choose{k+1}}\\
 =&
  \sum\limits_{k=0}^{a-2}
    \left(
     {{b+a-2-s}\choose{k+2}} -
     {{b+k-1-s}\choose{k+2}}
    \right)\\
 =&
  \sum\limits_{k=0}^{a-2}
   {{b+a-2-s}\choose{k+2}} -
  \sum\limits_{k=0}^{a-2}
   {{b+k-1-s}\choose{k+2}}
 \end{array}
}\]

Third, note that
\begin{equation}
{\allowdisplaybreaks
 \begin{array}{rl}
  \sum\limits_{k=0}^{a-2}{{b+a-2-s}\choose{k+2}} = &
  \sum\limits_{k=2}^{a}{{b+a-2-s}\choose{k}}\\
=&
  \sum\limits_{k=0}^{a-1}{{b+a-2-s}\choose{k}} +
  {{b+a-2-s}\choose{a}} - 1 - (b+a-2-s)\label{eq:Third1}\\
=&
  \sum\limits_{k=0}^{a-1}{{b+a-2-s}\choose{k}} +
  {{b+a-2-s}\choose{a}} - b - a + 1 + s
 \end{array}
}\end{equation}
and
\begin{equation}
{\allowdisplaybreaks
 \begin{array}{rl}
  \sum\limits_{k=0}^{a-2}{{b+k-1-s}\choose{k+2}} = &
  \sum\limits_{k=0}^{a-1}{{b+k-1-s}\choose{b-3-s}}
=
  \sum\limits_{p=b-1-s}^{b+a-3-s}{{p}\choose{b-3-s}}\\
\\
=&
 \begin{cases}
  {{b+a-2-s}\choose{b-2-s}}-{{b-1-s}\choose{b-2-s}}=
  {{b+a-2-s}\choose{a}} + s + 1 - b,\\
   \hspace{3cm} \text{if $(b-1-s\ge0)\iff(b>s)$}\\
  \ ~ \label{eq:Third2} %\tag{\label{dd}}\\
  \sum\limits_{p=0}^{b+a-3-s}{{p}\choose{b-3-s}} =
  {{b+a-2-s}\choose{b-2-s}} = {{b+a-2-s}\choose{a}},\\
   \hspace{3cm} \text{if $(b-1-s{}<0)\iff(b\le{}s)$}
 \end{cases}
 \end{array}
}\end{equation}

Finally, the equality for $s-a<b$ follows from
Equations~\eqref{eq:f(a,s,b)},~\eqref{eq:Third1} and \eqref{eq:Third2}
depending on either $s<b$ or $s\ge{}b$.

\emph{Case 1: $s<b$}:
\[
 \begin{array}{rl}
 N(a,s,b) =&
   a +\\
&
   \left( \sum\limits_{k=0}^{a-1} {{b+a-2-s}\choose{k}} +
          {{b+a-2-s}\choose{a}} - b - a + 1 + s
   \right) -\\
&
   \left( {{b+a-2-s}\choose{a}} + s + 1 - b
   \right)\\
=&
   \sum\limits_{k=0}^{a-1} {{b+a-2-s}\choose{k}}
 \end{array}
\]
\emph{Case 2: $s\ge{}b$:}
\[
 \begin{array}{rl}
 N(a,s,b) =&
   (a + b - 1 - s) +\\
&
   \left( \sum\limits_{k=0}^{a-1} {{b+a-2-s}\choose{k}} +
          {{b+a-2-s}\choose{a}} - b - a + 1 + s
   \right) -\\
&
   {{b+a-2-s}\choose{a}}\\
=&
   \sum\limits_{k=0}^{a-1} {{b+a-2-s}\choose{k}}.
 \end{array}
\]\hfill $\square$

\begin{definition}
\label{def:G(a,b)}
Let $G(a,b)$ be the number of matrices given by $N(a,s,b)$ such that
$0\le{}m_{1,1}\le{}a$, $0\le{}m_{r,2}\le{}b$,
$a\le{}m_{1,1}+m_{1,2}\le{}a+b$ (equivalently, $0\le{}m_{1,2}\le{}b$)
and, furthemore, $m_{1,1}>0$.
\end{definition}

\begin{lemma}
\label{lem:G(a,b)}
Let $G(a,b)$ be as Definition~\ref{def:G(a,b)}, then
\[
 G(a,b) = \sum\limits_{i=0}^{a}
           \sum\limits_{k=0}^{i-1}
            {{b}\choose{k+2}} + a\,b.
\]
\end{lemma}

\noindent \emph{Proof:} Note that in order to fulfill $m_{1,1}>0$ we
have to subtract $(a+b+1)$ matrices:
\[
\begin{array}{cll}
 (i\hspace{.25cm}b) &\quad\text{with $0<{}i\le{}a$}
                    &\text{\quad$\rightarrow$\quad$a$ matrices},\\
 (0\hspace{.25cm}j) &\quad\text{with $0<{}j\le{}b$}
                    &\text{\quad$\rightarrow$\quad$b$ matrices},\\
 (0\hspace{.25cm}0) &
                    &\text{\quad$\rightarrow$\quad$1$ matrix}.
\end{array}
\]

Now, using elemetal algebraic manipulations, the
equality~\cite{Cam94} follows:
\[
{\allowdisplaybreaks
\begin{array}{rl}
  G(a,b) =& \left(\sum\limits_{i=0}^a
                   \sum\limits_{j=0}^b
                    \sum\limits_{s=i}^{i+j} N(i,s,j)
            \right) - (a + b + 1)\\
         =& \left(\sum\limits_{i=0}^a
                   \sum\limits_{j=0}^b
                    \left(\left(\sum\limits_{s=i}^{i+j-1}
                                 \sum\limits_{k=0}^{i-1}{{j+i-2-s}\choose{k}}
                          \right) + 1
                    \right)
            \right) - (a + b + 1)\\
         =& \sum\limits_{i=0}^a
             \sum\limits_{j=0}^b
              \sum\limits_{k=0}^{i-1}
               \sum\limits_{s=i}^{i+j-2}
                {{j+i-2-s}\choose{k}} +
            \left(\sum\limits_{i=0}^a
                   \sum\limits_{j=0}^b
                    1
            \right) - (a + b + 1)\\
         =& \sum\limits_{i=0}^a
             \sum\limits_{j=0}^b
              \sum\limits_{k=0}^{i-1}
               \sum\limits_{s=i}^{i+j-2}
                {{j+i-2-s}\choose{k}} +
            (a+1)\,(b+1) - (a + b + 1)\\
         =& \sum\limits_{i=0}^a
             \sum\limits_{j=0}^b
              \sum\limits_{k=0}^{i-1}
               \sum\limits_{z=0}^{j-2}{z\choose{}k} +
            a\, b
         =  \sum\limits_{i=0}^a
             \sum\limits_{j=0}^b
              \sum\limits_{k=0}^{i-1}
               {j-1\choose{}k+1} +
            a\, b\\
         =& \sum\limits_{i=0}^a
             \sum\limits_{k=0}^{i-1}
              \sum\limits_{j=0}^{b-1}
               {j\choose{}k+1} +
            a\, b
         =  \sum\limits_{i=0}^a
             \sum\limits_{k=0}^{i-1}
              {b\choose{}k+2} +
            a\, b.
\end{array}
}\]\hfill $\square$

Finally, Lemma~\ref{lem:H(n)} proves Theorem~\ref{th:H(n)}.

\begin{lemma}
\label{lem:H(n)}
Let $H(n)$ and $G(a,b)$ as defined above, then
\[ H(n) = \sum\limits_{a=1}^n G(a, n-a) \]
and, furthermore,
\begin{eqnarray}
 H(n) = F(n+6) - (n^2 + 4n + 8).\notag% \label{eq:H(n)}
\end{eqnarray}
\end{lemma}

\noindent \emph{Proof:} From Definition~\ref{def:H(n)}
and~\ref{def:G(a,b)}, where $H(n)$ and $G(a,b)$ are respectively
defined, it is clear that
\[
 H(n) = \sum\limits_{a=1}^n G(a, n-a).
\]

Thus, we just have to prove that
\[ \sum\limits_{a=1}^n G(a, n-a) = F(n+6) - (n^2+4n+8). \]

In fact, taking into account that $G(0,n)=0$,
Lemma~\ref{lem:G(a,b)}, and the known identity of Fibonacci
numbers~\cite{Cam94}
\[ F(n) = \sum\limits_{k=0}^{n-1}{{n-k}\choose{k}}, \]
the equality results
\[
{\allowdisplaybreaks
\begin{array}{rl}
  \sum\limits_{a=1}^n G(a, n-a)
=&\sum\limits_{a=0}^n G(a, n-a)\\
=&\sum\limits_{a=0}^n
   \left(
    \sum\limits_{i=0}^a
     \sum\limits_{k=0}^{i-1}
      {n-a\choose{k+2}} +
    a\,(n-a)
   \right)\\
=&\sum\limits_{i=0}^n
   \sum\limits_{k=0}^{i-1}
    \sum\limits_{y=0}^{n-i}
     {{y}\choose{k+2}} +
  {{n+1}\choose{3}}
= \sum\limits_{i=0}^n
   \sum\limits_{k=0}^{i-1}
    {{n-i+1}\choose{k+3}} +
  {{n+1}\choose{3}}\\
=& \sum\limits_{k=0}^{n-1}
    \sum\limits_{i=k+1}^{n}
     {{n-i+1}\choose{k+3}} +
   {{n+1}\choose{3}}
= \sum\limits_{k=0}^{n-1}
   \sum\limits_{x=1}^{n-k}
    {{x}\choose{k+3}} +
  {{n+1}\choose{3}}\\
=&\sum\limits_{k=0}^{n-1}
   {{n-k+1}\choose{k+4}} +
  {{n+1}\choose{3}}
= \sum\limits_{k=4}^{n+3}
   {{n+5-k}\choose{k}} +
  {{n+1}\choose{3}}\\
=&\sum\limits_{k=0}^{n+5}
   {{n+5-k}\choose{k}} -
  \left[{{n+5}\choose{0}} +
        {{n+4}\choose{1}} +
        {{n+3}\choose{2}} +
        {{n+2}\choose{3}}
  \right] +
  {{n+1}\choose{3}}\\
=&F(n+6) - (n^2 + 4n + 8).
\end{array}
}\]\hfill $\square$

\section*{Acknowledgements}
Josep Freixas was partially supported by Grant MTM 2006--06064
of ``Ministerio de Educaci\'on y Ciencia y el Fondo Europeo de
Desarrollo Regional'', SGRC 2009-1029 of ``Generalitat de
Catalunya".

Xavier Molinero and Salvador Roura were partially supported by
the Spanish ``Ministerio de Educaci\'on y
Ciencia'' programme TIN2006-11345 (ALINEX) and SGRC 2009-1137 of
``Generalitat de Catalunya".

Josep Freixas, Xavier Molinero and Salvador Roura were also
partially supported by Grant 9-INCREC-11 of
``(PRE/756/2008) Ajuts a la iniciaci\'o/reincorporaci\'o a la recerca
(Universitat Polit\`ecnica de Catalunya)''.

\bibliographystyle{plain}
%\bibliography{references}

\end{document}